\documentclass[a4paper,12pt]{amsart}
\usepackage{amsmath}
\usepackage{graphicx}
\usepackage[latin1]{inputenc}
\usepackage{epstopdf}
\usepackage{amssymb,amscd,array}
\usepackage{color}
\usepackage{float}

\newtheorem{thm}{Theorem}[section]
\newtheorem{cor}[thm]{Corollary}
\newtheorem{lem}[thm]{Lemma}
\newtheorem{prop}[thm]{Proposition}
\theoremstyle{definition}
\newtheorem{defn}[thm]{Definition}
\theoremstyle{remark}

\numberwithin{equation}{section}
\newcommand{\eps}{\varepsilon}

\makeindex \makeatletter
\def\captionof#1#2{{\def\@captype{#1}#2}}
\makeatother

\newcounter{tablegroup}
\newcounter{subtable}[tablegroup]

\begin{document}
\title{On totally periodic $\omega$-limit sets}

\author{ Habib Marzougui and Issam Naghmouchi}

\address{ Habib Marzougui, University of Carthage, Faculty
of Sciences of Bizerte, Department of Mathematics,
Jarzouna, 7021, Tunisia; Issam Naghmouchi, University of Carthage, Faculty
of Sciences of Bizerte, Department of Mathematics,
Jarzouna, 7021, Tunisia.}
 \email{hmarzoug@ictp.it and habib.marzougui@fsb.rnu.tn; issam.nagh@gmail.com and issam.naghmouchi@fsb.rnu.tn}
 \thanks{This work is supported by the research laboratory: syst\`emes dynamiques et combinatoire: 99UR15-15}


\begin{abstract}
An $\omega$-limit set of a continuous self-mapping of a compact metric space $X$ is said to be totally periodic if all of
its points are periodic. We say that $X$ has the $\omega$-FTP property provided that for each continuous self-mapping $f$ of $X$,
every totally periodic $\omega$-limit set is finite. Firstly, we show that connected components of every totally periodic $\omega$-limit
set are
finite. Secondly, we show in one hand, that a zero-dimensional compact metric space has the $\omega$-FTP property, and in the other hand,
for the wide class of one-dimensional continua, we prove that
a hereditary locally connected $X$ has the
$\omega$-FTP property if and only if $X$ is completely regular. This holds in particular for $X$ being
 a local dendrite with discrete set of branch points, and in particular, for a graph. For higher dimension, we show that
any compact metric space $X$ containing a free topological $n$-ball $(n\geq 2)$ does not admit the $\omega$-FTP property. This holds
in particular, for any topological compact manifold of dimension greater than $1$.

\end{abstract}
\maketitle

\section{\bf Introduction}

 Let  $(X,d)$ be a compact metric space and $f$ be a self continuous map of $X$. Let $\mathbb{Z}_{+}$ and $\mathbb{N}$ be the
sets of non-negative integers and positive integers respectively.
Denote by $f^{n}$ the $n$-th iterate of $f$; that is, $f^{0} = \textrm{id}_{X}$: the identity and $f^{n} = f\circ f^{n-1}$ if
$n\in\mathbb{N}$. For any $x\in X$ the subset $O_{f}(x) = \{f^{n}(x):
\ n\in\mathbb{Z}_{+}\}$ is called the $f$-orbit of $x$. A point
$x\in X$ is called periodic of prime period $n\in\mathbb{N}$ if
$f^{n}(x)=x$ and $f^{i}(x)\neq x$ for $1\leq i\leq n-1$, the orbit
of such point is called periodic orbit. We denote by P($f)$ the set
of periodic points and by Fix$(f)$ the set of fixed points of $f$. A
subset $A\subset X$ is called \textit{$f$-invariant} if $f(A)\subset
A$, it is strongly $f$-invariant if $f(A)=A$. A $p$-tuple $(A_{0}, \dots, A_{p-1})$ of subsets of $X$ is called a \textit{periodic cycle} if
$f(A_{0})=A_{1}, f(A_{1})=A_{2}, \dots, f(A_{p-1})=A_{0}$. For a subset $A$ of
$X$, denote by $\overline{A}$ the closure of $A$. We define the $\omega$-limit set of a point $x\in X$ to be the set
$\omega_{f}(x) = \{y\in X: \exists\ n_{i}\in \mathbb{N},
n_{i}\to +\infty, \underset{i\to +\infty}\lim d(f^{n_{i}}(x), y)
= 0\}$.
\medskip

 The topological characterization of $\omega$-limit sets of $f$ is known only in few spaces, such as zero-dimensional spaces, this was found in
\cite{sp}, Theorem 13), graphs (which includes the intervals, circles), this was given by \cite{hr} and more generally,
for hereditally locally connected continuum (see below for the definition),
this was given by \cite{sp}. For higher-dimensional spaces, only
partial results are known; among them $\omega$-limit sets in the square \cite{ag}, \cite{jl}, \cite{ba}.
Apart from these results, the topological characterization of $\omega$-limit sets is unknown. In this note, we focus on the class of $\omega$-limit sets that are composed of periodic points.
We say that $\omega_{f}(x)$ is \textit{totally periodic} if \ $\omega_{f}(x)\subset P(f)$.
The topological characterization of such sets is not yet been examined.
The question whether a totally periodic $\omega$-limit set is finite, is addressed.
\medskip

\begin{defn}\label{def2}
 A compact metric space $X$ has the $\omega$-FTP property provided that for each continuous self-mapping $f$ of $X$, every totally periodic
$\omega$-limit set is finite (i.e. a periodic orbit).
 \end{defn}
\medskip

 The first result show that every totally periodic $\omega$-limit set has finitely many connected components.
\medskip

 \begin{thm}\label{thm1} Let $X$ be a compact metric space, $f$ be a continuous self-mapping of $X$ and $x\in X$. If $\omega_{f}(x)$ is totally periodic then
$\omega_{f}(x)$ has finitely many connected components that form a periodic cycle.
\end{thm}
\medskip

\begin{cor}\label{n3}
Let $X$ be a compact metric space, $f$ be a self continuous map of $X$. Let $n\geq 1$ be an integer and $x\in X$. If
$\omega_{f}(x)\subset \mathrm{Fix}(f^n)$ then the number of connected components of $\omega_{f}(x)$ divides $n$.
\end{cor}
\medskip

Recall that a continuum is a compact connected metric space. We give a complete answer of the question: whether a continuum has the $\omega$-FTP property, in the wide class of one-dimensional continua : the  hereditary locally connected continua.
\medskip

A space is called \textit{degenerate} provided it has only one point; otherwise it is called \textit{non-degenerate}. A continuum $X$ is said to be:
\medskip

 \textbullet \ hereditarily locally connected if every subcontinuum of $X$ is locally connected;

%\textbullet \ regular if $X$ has a basis of open sets with finite boundaries;

\textbullet \ completely regular if every non-degenerate subcontinuum of $X$ has non-empty interior.

\textbullet  \ dendrite if $X$ is locally connected and contains no simple closed curve;

\textbullet  \ local dendrite if every point of $X$ has a neighborhood which is a dendrite,

\textbullet  \ graph if $X$ can be written as the union of finitely many arcs any two of which are either disjoint or intersect only in
one or both of their end points.
\smallskip

Recall that if $X$ is a local dendrite, a point $x\in X$ is called a \textit{branch point} (\textit{resp.} an \textit{end point}) if for some neighborhood $U$ of $x$, $U$ is a dendrite and $U\backslash \{x\}$ has more than two connected components (\textit{resp.} $U\backslash \{x\}$ is connected).
Denote by $E(X)$ and $B(X)$ the set of all end points and branch points of $X$, respectively.
It is known that every completely regular continuum is
hereditarily locally connected (\cite{ku}, Theorem 50.IV.1). Note that any graph as well as any dendrite is hereditarily locally connected.
Thus every subcontinuum of a graph (resp. a dendrite) is a graph (resp. a dendrite).
Any graph as well as any local dendrite with branch points discrete
is completely regular. More information about hereditarily locally connected continua can be found in \cite{Nadler}.
Our second result can be stated as follows.
\medskip

 \begin{thm}\label{thm2} Let $X$ be hereditally locally connected continuum. Then $X$ has the $\omega$-FTP property if and only if $X$
is completely regular.
\end{thm}
\medskip

In particular:
\medskip

\begin{cor}\label{c2}
If $X$ is a local dendrite with $B(X)$ discrete then $X$ has the $\omega$-FTP property. In particular, this holds whenever $X$ being either:

\textbullet \ a graph,

or

\textbullet \ a dendrite with $B(X)$ discrete (in particular, a dendrite with $E(X)$ closed).
\end{cor}
\medskip

 \textbf{Remark 1}.  If $X$ is a dendrite and $f$ is a monotone continuous
map of $X$ into itself, then any totally periodic $\omega$-limit set is finite. This follows from (\cite{Nagh1}, Theorem C).
 \bigskip

 \textbf{Remark 2}. In \cite{Li}, Li introduced the definition of $\omega$-scrambled set for a continuous map $f$ of a compact metric space $X$ into itself 
 as follows:
A subset $S$ of $X$ is called $\omega$-scrambled for $f$ if for any $x, \ y\in S$ with $x\neq y$:
\begin{itemize}
 \item[(i)] $\omega_{f}(x) \backslash \omega_{f}(y)$ is uncountable,

\item[(ii)]  $\omega_{f}(x)\cap \omega_{f}(y)$ is nonempty,

\item[(iii)] $\omega_{f}(x) \backslash P(f)$ is nonempty.
\end{itemize}
\medskip

Note that condition (iii) is equivalent to say that $\omega_{f}(x)$ is not totally periodic. 
The definition of $\omega$-scrambled set is reduced only to conditions (i) and (ii) 
when the space $X$ is either a completely regular continuum or a zero-dimensional compact space; this results from
Theorem \ref{thm2} and Corollary \ref{c1}.

\medskip

In higher dimension, let $n\geq 2$ be an integer and let $B_n$ denote the unit $n$-ball given by $B_n = \{x\in\mathbb{R}^n: \ \|x\|\leq 1\}$,
where $ \| \ \|$ in the Euclidean norm on $\mathbb{R}^n$, $S^{n-1}=\partial B_{n}$ its boundary and $\textrm{int}(B_{n})=B_n\backslash \partial B_{n}$ its interior. A topological $n$-ball is a topological
space homeomorphic to $B_n$.
Let $X$ be a topological space. A topological $n$-ball $B\subset X$ is called \textit{free} in $X$ if $h(\textrm{int}(B_{n}))$
is open in $X$ where $h$ is any homeomorphism from $B_n$ onto $B$.
Our third main result is the following.
\medskip

\begin{thm}\label{t4} Let $X$ be a compact metric space containing a free topological $n$-ball $(n\geq 2)$. Then there exists a homeomorphism of $X$ into itself having an infinite
$\omega$-limit set consisting of fixed points. In particular, $X$ does not admit the $\omega$-FTP property.
\end{thm}
\medskip

\begin{cor}\label{c18}
Every topological compact manifold of dimension $\geq 2$ does not admit the $\omega$-FTP property.
\end{cor}
\medskip

However, the following holds:
\medskip

\begin{prop}\label{p3} Let $X_{1},\dots, X_{n}$ are completely regular continua and let $f$ be the self map of
$X_{1}\times \dots \times X_{n}$ defined by
$f((x_{1}, \dots, x_{n})) = (f_{1}(x_{1}), \dots,f_{n}(x_{n}))$ where $f_{1}, \dots f_{n}$ are continuous self-mapping of $X_{i}$.
Then every totally periodic $\omega$-limit set of $f$ is finite.
 \end{prop}
\medskip

The proposition \ref{p3} holds, for example, whenever $X_{i}=[0,1]$ or $X_{i}= S^{1}$.
\medskip

This paper is organized as follows. In Section 2 we prove Theorem \ref{thm1} and
Corollary \ref{n3}. Section 3 is devoted to the one dimensional case, in particular we prove Theorem \ref{thm2}. Finally,
in Section 4, we study the higher dimension by proving Theorem \ref{t4} and Proposition \ref{p3}.
\medskip

%We briefly recall some results which will be needed in the rest of the paper.

\section{\bf The connected component of a totally periodic  $\omega$-limit set}
\medskip

Let $X$ be a compact metric space, $f$ be a self continuous map of $X$ and $x\in X$.
We recall that the $\omega$-limit sets possess the following basic properties:
\medskip

\begin{prop}[\cite{lsB}, Lemma 2, Chapter IV]\label{p1}
The set $\omega_{f}(x)$ is a non-empty, closed and strongly
invariant set.
\end{prop}
\medskip

\begin{prop}[\cite{lsB}, Lemma 4, Chapter IV]
An $\omega$-limit set $\omega_{f}(x)$ is finite if and only if it is
the orbit of some periodic point.
\end{prop}
\medskip

\begin{lem}[\cite{lsB}, Lemma 3, p. $71$] \label{n1}  Set
$L=\omega_{f}(x)$ and $F$ be any non-empty
 proper closed subset of $L$. Then $F \cap
 \overline{f(L\backslash F)} \neq \emptyset$.
\end{lem}
\medskip

\begin{lem}[\cite{K}, Proposition 3.5, p. $107$] \label{n2}
If $X$ is a compact metric space and totally disconnected then the clopen sets form a base for its topology.
\end{lem}
\medskip

\begin{lem} \label{prop1}
Assume that $\omega_{f}(x)$ is totally periodic. If $\omega_{f}(x)$ is totally disconnected then it is finite.
\end{lem}
\medskip

\begin{proof}
 Write $L = \omega_{f}(x)$. By Proposition \ref{p1}, we have $f(L)=L$. We claim that
the restriction map $f_{| L}:L \rightarrow L$ is an
homeomorphism:
 Let $x,y\in L$ with period $p$ and $q$ respectively.
If $f(x)=f(y)$ then $f^{pq}(x) = f^{pq}(y)$ hence $x = y$ then
$f_{|L}$ is bijective. As $L$ is compact then $f_{|L}$ is a homeomorphism.
By hypothesis, $L\subset P(f)$, so $ L =  \bigcup\limits_{n=1}^{+\infty}
\textrm{Fix}(f^{n})\cap L$. Hence $\textrm{Fix}(f^{p})
\cap L$ has non-empty interior in $L$ for some $p\geq 1$. By Lemma
\ref{n2}, there exists a non-empty clopen subset $U$ of $L$
such that $U \subset \textrm{Fix}(f^{p})$. Denote by $F=U \cup f(U)
\cup \dots \cup f^{p-1}(U)$. Then $F$ is clopen in $L$ and we have $f(F)= F$.

 Suppose that $L$ is infinite. One can choose $U$ so that $F\neq L$: indeed, if $U=\{y\}$, for some $y\in U$ then $F= \{y,\dots, f^{p-1}(y)\}$
and $F\neq L$. If  $U\neq \{y\}$ for some $y\in U$, we choose a non-empty clopen subset $U^{\prime}\subset U\backslash\{y\}$ and so
 $F:=U^{\prime} \cup f(U^{\prime})
\cup \dots \cup f^{p-1}(U^{\prime})\varsubsetneq U \cup f(U)
\cup \dots \cup f^{p-1}(U)$, hence  $F\neq L$.

Now by Lemma \ref{n1}, we have
$\overline{f(L\backslash F)} \cap F=f(L\backslash F) \cap F
\neq \emptyset $. Hence, there is a point $z\in L\backslash F$ such that $f(z)\in F$. As $z\in P(f)$ and  $f(F)= F$, thus $z\in F$, a contradiction. We conclude that $L$ is finite.
\end{proof}
\medskip

 As a consequence, we have the following result.

\begin{cor}\label{c1}
If $X$ is a compact metric space and totally disconnected, then $X$ has the $\omega$-FTP property. In particular, Cantor space has the $\omega$-FTP property.
\end{cor}
\medskip

 Now let $X$ be a metric compact space, $f$ be a self continuous map of $X$ and assume that $L=\omega_{f}(x)$ is totally periodic.
Denote by $Y=\overline{O_{f}(x)}=O_{f}(x) \cup L$. We have $f(Y)\subset Y$. For every $z\in Y$, $[z]$ means the connected component in $Y$ that
 contains $z$. Let $\mathcal{C}$ be the family of connected components of $Y$ which form a decomposition of $Y$.
  By collapsing these components to points, we get the quotient space $Y/\mathcal{C}$ which is a compact metric
  space and totally disconnected (cf. \cite{bdr}).
 We denote by $\pi: Y\longrightarrow Y/\mathcal{C}; y\longmapsto [y]$ and
 $\displaystyle{\tilde{f}:Y/\mathcal{C}\rightarrow
  Y/\mathcal{C}}$ the induced map defined by: $\forall \ [y]\in Y/\mathcal{C}, \quad \tilde{f}([y])=[f(y)]$ i.e.
$\tilde{f}\circ \pi = \pi \circ f$. It is plain that $\tilde{f}$ is well defined and continuous.
 \medskip

\begin{lem}\label{l33}
The following assertions hold:

 \begin{enumerate}
 %\item The map $f_{|L}:L \rightarrow L$ is a homeomorphism.
 \item for any $n\in \mathbb{Z_+}$,  $[f^{n}(x)] = \{f^{n}(x)\}$ and $[f^{n}(x)]$ is an isolated point in $Y/\mathcal{C}$.

\item If $y\in L$ then $[y]\subset L$.
 \item for any $y\in Y$, $f([y]) = [f(y)]$ and $\forall n\in \mathbb{Z_+}, f^{n}([y])=[f^{n}(y)]$.
 \item for any $z\in L$, $[z]$ is a periodic point for $\tilde{f}$.

\item If $L$ is infinite then $\{[y], y\in L\} = \omega_{\tilde{f}}([x])$.
 \end{enumerate}
\end{lem}
\medskip

\begin{proof} (1) If $f^{n}(x)\in L$ then $f^{n}(x)\in P(f)$ and so $Y = O_{f}(x)$ is finite. Therefore $f^{n}(x)$ is isolated in $Y$.
If $f^{n}(x)\notin L$, there is $\varepsilon_{1}>0$ such that
$B(f^{n}(x),\varepsilon_{1}) \cap L = \emptyset
 $, where $B(x_{0}, r)$ denotes the open ball of radius $r>0$ and center $x_{0}\in X$. In the other hand,
 there is $\varepsilon_{2} >0$ such that $B(f^{n}(x),\varepsilon_{2}) \cap \mathcal{O}_{f}(x)=\{f^{n}(x)\}$
(otherwise $f^{n}(x) \in L$, a contradiction). So $B(f^{n}(x),\varepsilon) \cap Y=\{f^{n}(x)\} $ where
$\varepsilon= \min(\varepsilon_{1},\varepsilon_{2})$. Thus $f^{n}(x)$ is isolated in $Y$ and hence $
 [f^{n}(x)]=\{f^{n}(x)\}$. As
$\pi^{-1}(\{[f^{n}(x)]\})=[f^{n}(x)] = \{f^{n}(x)\}$ which is open in $Y$, it follows that $[f^{n}(x)]$
 is isolated in $Y/\mathcal{C}$.
\\\\
(2) Let $y\in L$. If for some $n\in\mathbb{Z_+}$, $f^n(x)\in [y]$ then by assertion (1), $[y]=[f^n(x)]=\{f^n(x)\}$. 
So $y=f^n(x)$ and hence $[y]=\{y\}\subset L$. If now $O_f(x)\cap [y]=\emptyset$ then it is clear that $[y]\subset L$.
\\\\
(3) By assertion (1), if $y\in O_f(x)$ then it is clear that $f([y])=[f(y)]$. Now let $y\in L$. Since $[y]$ is the connected component of $Y$ containing $y$, $f$ is continuous and  $f(Y)\subset Y$, so $f([y])\subset Y$ is
connected and contains $f(y)$, therefore, $f([y])\subset [f(y)]$ and so $f^{i}([y])\subset [f^{i}(y)]$, for every $i\in 
\mathbb{N}$.
 Conversely, suppose that $f([y]) \varsubsetneq [f(y)]$. 
By assertion (2) and as $f(L)= L$, 
$[f(y)]\subset L$ and there exists $z\in L \backslash [y]$ such that $f(z)\in [f(y)]$. Let $p, \ k \geq 1$
such that $f^{p}(y)= y$ and $f^{k}(z)= z$. Then we get $f^{kp}(z)=z \in [f^{kp}(y)] = [y]$, a contradiction.
\\\\
(4) If $z\in L$ with period $k\geq 1$ then $\tilde{f}^{k}([z])= [f^{k}(z)]=[z]$ and so $[z]$ is a periodic point of $\tilde{f}$.
\\\\
(5) Since $L= \omega_{f}(x) = \underset{n\in \mathbb{N}}\cap\overline{\{f^{k}(x): k\geq n\}}$, we have $\pi(L)= \{[y],  y \in L\} \subset \omega_{\tilde{f}}([x])$. Conversely, let
$[y]\in \omega_{\tilde{f}}([x])$. If $y=f^{k}(x)$ for some $k\in
 \mathbb{Z}_+$, then $\omega_{\tilde{f}}([x])= \omega_{\tilde{f}}([y])$. Now by ($1$), $[y]$ is isolated in
$Y/\mathcal{C}$ and  $\{f^k(x)\}=[y]$. Hence $[y]= \tilde{f}^{k+p}([x])$
for infinitely many $p$. Therefore $[y]$ is $\tilde{f}$-periodic and so $f^k(x)$ is $f$-periodic point. Thus $L$ is finite, a contradiction.
We conclude that $y\in L$ and therefore $\omega_{\tilde{f}}([x])=
\{[y], \ y \in L\}$.
\end{proof}
\bigskip

\textit{Proof of Theorem \ref{thm1}}.  If $L$ is finite, the proof is clear. Assume that  $L$ is infinite. As $Y/\mathcal{C}$ is compact and totally disconnected metric space, and  $\omega_{f}(x) \subset P(f)$, then
 $\omega_{\tilde{f}}([x])$ is totally disconnected and by Lemma \ref{l33}, (4) and (5), $\omega_{\tilde{f}}([x]) \subset P(\tilde{f})$.  Therefore by Lemma \ref{prop1}, $\omega_{\tilde{f}}([x])$ is
finite i.e. $\omega_{\tilde{f}}([x])= \{[z], \tilde{f}([z]),\dots,\tilde{f}^{k-1}([z])\}$ (for some $z \in Y$ and $k\geq 1$). As
$\tilde{f}^{i}([z])= f^{i}([z])$ (Lemma \ref{l33}, (3)) we conclude that $L$ has a
 finite number of connected components: $[z],f([z]),\dots, f^{k-1}([z])$ that form a periodic cycle for $f$.
The proof is complete. \qed
\bigskip

%\textit{Proof of Corollary \ref{c1}}. By Theorem \ref{thm1}, $L$ has a finite number of connected components, then $L$ is finite
%(i.e. a periodic orbit) since $X$ is totally disconnected. \qed
%\bigskip
\smallskip

\textit{Proof of Corollary \ref{n3}}.
As $\omega_f(x)=\omega_{f^n}(x)\cup\omega_{f^n}(f(x))\dots\cup \omega_{f^n}(f^{n-1}(x))$ then for each $0\leq i\leq n-1$,
$\omega_{f^n}(f^i(x))\subset \mathrm{Fix}(f^n)$. By Theorem \ref{thm1},
$\omega_{f^n}(f^i(x))$ has finitely many connected components $C_{1}, \dots, C_{p}$, so
each $C_{k}$, ($1\leq k\leq p$) is open in $\omega_{f^n}(f^i(x))$. Suppose that $\omega_{f^n}(f^i(x))$ is not connected, then $p\geq 2$ and by applying
 Lemma \ref{n1}, there exists $z\in \omega_{f^n}(f^i(x))\backslash C_{1}$ such that $f^{n}(z)\in C_{1}$, a contradiction since  $f^{n}(z)=z$.
We conclude that for each $1\leq i\leq n$, $\omega_{f^n}(f^i(x))$ is connected and therefore the number $l$ of connected components of
$\omega_f(x)$ is at most $n$. As the connected components of
$\omega_f(x)$ form a periodic cycle for $f$ then $l$ divides $n$. \qed
\bigskip

\section{\bf The $\omega$-FTP property in one-dimension}

\textbf{3.1. The example.} We are going to construct a continuous map $F_{0}$ on a dendrite $D_{0}$ admitting an $\omega$-limit set
 which is an arc formed only by fixed points.

\textbf{(a)} \ First, we define the space $D_{0}$ as a subset of the plane in the following way:

Let $n\in \mathbb{N}$. Denote by $S_{n} = \{\frac{i}{2^n}: \ 1\leq i\leq 2^n-1 \}$. For $n=1$, we let $\Omega_1=\{\frac{1}{2}\}$
and for each $n>1$, we let $$\Omega_{n} = S_{n} \backslash \cup_{i=1}^{n-1}\Omega_{i}.$$
For example, $\Omega_2 = \{\frac{1}{4}, \frac{3}{4}\}$, $\Omega_3 = \{\frac{1}{8}, \frac{3}{8}, \frac{5}{8}, \frac{7}{8}\}$.
For each $n\in\mathbb{N}$, the cardinal of $\Omega_n$ is $2^{n-1}$.
We define now a dense sequence $(a_n)_{n\in\mathbb{N}}$ in $[0,1]$ as follow:

 For $n=1$, we let $a_1=\frac{1}{2}$. For each $n\in\mathbb{N}$, set

\textbullet \ $\Omega_{2n}= \{a_i: 2^{2n-1}\leq i\leq 2^{2n}-1\}$ given in this order: \ $a_{2^{2n-1}}<\dots< a_{2^{2n}-1}$.

\textbullet \ $\Omega_{2n+1}=\{a_i: 2^{2n}\leq i\leq 2^{2n+1}-1\}$  given in this order: \ $a_{2^{2n}}>\dots> a_{2^{2n+1}-1}$.

For example, $a_{2}= \frac{1}{4}$,  $a_{3}= \frac{3}{4}$, $a_{4}= \frac{7}{8}$, $a_{5}= \frac{5}{8}$, $a_{6}= \frac{3}{8}$, $a_{7}= \frac{1}{8}$.

For any $n, k\in\mathbb{N}$ such that $a_{k}\in \Omega_n$, we let $I_k = \{a_k\}\times [0,\frac{1}{n}]$. The dendrite $D_{0}$ is then given by:
 $$D_{0}:=([0,1]\times \{0\})\cup \bigcup_{k\in\mathbb{N}}I_k.$$
%(see Figure 1).
\medskip

\textbf{(b)} \ Second we construct the map $F_{0}$ as follows: For any $k\in\mathbb{N}$ such that $a_{k}\in \Omega_n$, write
$I_{k} = \{a_{k}\}\times [0,\frac{1}{n}] = [A_k, B_k]$,
where $A_k = (a_k,0)$,  $B_k = (a_k,\frac{1}{n})$. We let $C_k = (a_k,\frac{2}{n})$. Then the map $F_{0}$ is defined inductively as follows:
 \begin{itemize}
 \item $[0,1]\times\{0\}= \textrm{Fix}(F_{0})$.\\
\item $(F_{0})_{| [A_k,C_k]}: [A_k,C_k]\to [A_k,A_{k+1}]$ be a homeomorphism that fixes the point $A_k$ and maps linearly
$[A_k,C_k]$ to $[A_k,A_{k+1}]$,\\
\item $(F_{0})_{\mid [C_k,B_k]}:[C_k,B_k]\to [A_{k+1},B_{k+1}]$ be a homeomorphism that maps linearly $[C_k,B_k]$ to $[A_{k+1},B_{k+1}]$
and sends $C_k$ to $A_{k+1}$ and $B_k$ to $B_{k+1}$.
 \end{itemize}

The map $F_{0}$ is continuous on $D_{0}$: Let $k\in\mathbb{N}$ and $M\in I_k\setminus A_k$. Then there is a sub-arc of $I_k$
containing $M$, which is open in $D_{0}$. Hence $F_{0}$ is continuous at $M$ since $F_{0}$ is a homeomorphism on $I_k$. Let
 $M=(x,0)\in [0,1]\times \{0\}$ and $(M_n)$ be a sequence of points in $D_0$ that converges to $M$. If $M_n$ lies eventually in
 $[0,1]\times \{0\}$ then $ F_0(M_n)= M_{n}$ and so $\underset{n\to +\infty}\lim F_0(M_n)=M$. If $M_n$ lies eventually in
 $I_{k(n)}\setminus A_{k(n)}$ then we distinguish three subcases:\\
   \begin{itemize}
     \item Case 1. $M_n$ lies finitely many times in each $I_{k(n)}\setminus A_{k(n)}$. In this case,
$\underset{n\to +\infty}\lim k(n)=+\infty$. Since $\underset{n\to +\infty}\lim \textrm{diam}(I_{k(n)})=0$, we obtain
 $\underset{n\to +\infty}\lim A_{k(n)}=M$. We also have $\underset{n\to +\infty}\lim \textrm{diam} ([A_{k(n)},A_{k(n)+1}])=0$ and for
eventually $n$, $F_0(M_n)\in [A_{k(n)},A_{k(n)+1}]\cup I_{k(n)+1}$. Therefore $\underset{n\to +\infty}\lim F_0(M_n)=M$.

\item Case 2. $M_n$ lies in a finite union of $I_k$ (i.e. $k(n)$ is bounded). In this case,
 then there exists $k_0\in\mathbb{N}$ such that $M_n$ lies eventually in $I_{k_0}$ and $A_{k_0}=M$. By the continuity of the restricted map
$(F_{0})_{\mid I_{k_0}}$, we get $\underset{n\to +\infty} \lim F_0(M_n)=M$.

     \item Case 3. There exists $k_0\in\mathbb{N}$ such that $M_n$ lies infinitely many times in $I_{k_0}$ and
$M_n\notin I_{k_0}$ for infinitely many times. In this case, the sequence $(M_n)_{n}$ has two subsequences $(M_{\alpha(n)})_{n}$ and
$(M_{\beta(n)})_{n}$ such that  $\mathbb{N}=\{\alpha(n): \ n\in\mathbb{N}\}\cup \{\beta(n): \ n\in\mathbb{N}\}$ where $(M_\alpha(n))$ satisfies
case 1 and  $(M_{\beta(n)})$ satisfies the case 2. So $\underset{n\to +\infty}\lim F_0(M_n)=M$.
   \end{itemize}

In result, $F_0$ is continuous on $[0,1]\times \{0\}$ and so $F_0$ is continuous on the whole dendrite $D_0$. The $F_{0}$-orbit of
the point $B_1$ is then $O_{F_{0}}(B_{1})= \{B_k: \ k\in\mathbb{N}\}$. Hence
$\omega_{F_{0}}(B_1) = [0,1]\times\{0\} = \textrm{Fix}(F_{0})$.
\bigskip
\bigskip

%\begin{figure}[!h]
%\begin{center}
%\includegraphics[width=0.8\textwidth]{D0.jpg} \caption{Dendrite $D_0$}
%\end{center}
%\end{figure}
\bigskip
\bigskip
\bigskip
\bigskip
\bigskip
\bigskip
\textbf{3.2. Proof of Theorem \ref{thm2} }
\medskip

\textbullet \ Assume that $X$ is completely regular and let $L=\omega_{f}(x)\subset P(f)$. We claim that
$L$ is totally disconnected. Indeed, otherwise $L$ contains
a non-degenerate connected component $C$ necessarily with
non-empty interior in $X$. In particular $L$ is infinite and there exists $n\in\mathbb{N}$ such that $f^n(x)\in L$.
Thus $f^n(x)$ is a periodic point. As
$\omega_{f}(x) = \omega_{f}(f^n(x))$, then $L$ is finite, a
contradiction. Now by Lemma \ref{prop1}, $L$ itself is
finite.
\\\\
 \textbullet \ Conversely, suppose that $X$ is not completely regular.
\medskip

\textbf{(a)} First we prove that $X$ contains a homeomorphic copy of the dendrite $D_{0}$ (see the Example 3.1 above). Indeed, there is a non-degenerate locally connected
subcontinuum $Y$ of empty interior. As $Y$ itself is arcwise connected, it contains an arc $[a,b]$. Take $c\in (a,b)$ and $\varepsilon>0$.
As $X$ is locally arcwise connected, there is an arcwise connected neighborhood $U$ of $c$ in $X$ with diameter less than $\varepsilon$
 (see \cite{Nadler}, Theorem 8.25, p.131). As $Y$ has empty interior in $X$, there exists
$d\in U\setminus [a,b]$ and so an arc, say $[c,d]$ in $U$, joining $c$ and $d$. Let $r\in [c,d[$ be such that the sub-arc
$[c,r]$ of $[c,d]$ intersects $[a,b]$ only in the point $r$. Hence we have proved that for each $\eps>0$ and $c\in (a,b)$, there is $r\in [a,b]$ with $d(r,c)<\eps$
and an arc $I$ with $r$ as one of its endpoints such that $I\cap [a,b]=\{r\}$. So if we take two distinct points $c,c^{\prime}\in (a,b)$,
we can find two disjoint arcs one has $c$ as one of its endpoints and the other has $c^{\prime}$ as one of its endpoints and each one of them
intersects $[a,b]$ only in $c$ and $c^{\prime}$ respectively. We conclude the existence of a dense sequence $(r_n)_{n}$ of points in $[a,b]$
and a pairwise disjoint arcs $(I_n)_{n}$, each of them has $r_n$ as one of its endpoints with $I_{n}\cap [a,b]=\{r_n\}$.
Since $X$ is hereditarily locally connected, then by (\cite{ku}, Theorem 50.IV.9), the arcs $(I_n)_{n}$ is a null family; i.e. for any
$\varepsilon >0$, only finitely many $I_{n}$ have diameter $> \varepsilon$. Therefore $\underset{n\to +\infty}\lim \textrm{diam}(I_{n}) = 0$.
We conclude that the set $Z:= [a,b]\cup \underset{n\in\mathbb{N}}\cup I_n$ is homeomorphic to the dendrite $D_{0}$.
\medskip
\smallskip

\textbf{(b)} Second, we need the following proposition.

\begin{prop}[\cite{sp}, proposition 11]\label{p:4}
 Let $X$ be a hereditarily locally connected continuum and let $Y$ be a locally connected continuum. Then any
continuous map $f: M \longrightarrow Y$ defined on a closed subset $M$ of $X$ can be extended to a continuous map $F : X\longrightarrow Y$
defined on $X$.
\end{prop}
\medskip

Now, we complete the proof of Theorem \ref{thm2} as follows. From the Example 3.1, there exist a point $x\in Z$ and a
continuous map $F: Z\longrightarrow Z$ defined on $Z$ such that $\omega_F(x)= [a,b]= \textrm{Fix} F$.
So, by Proposition \ref{p:4}, the map $F$ can be extended to a continuous map $\tilde{F}$ defined on $X$ into $X$
with $\omega_{\tilde{F}}(x)=
[a,b]\subset \textrm{Fix}\tilde{F}$. Therefore $X$ does not admit
the $\omega$-FTP property. \qed
\medskip

\section{\bf The $\omega$-FTP property in higher dimension}
\medskip

\begin{prop}\label{prop51} For each integer $n\geq 2$, there exists a homeomorphism $f$ of $B_n$ onto itself such that
$S^{n-1} \cup \{0\}\subset \mathrm{Fix}f$ and $\omega_f(x)$ is a circle included in
$S^{n-1}$, for some $x\in B_n\setminus (S^{n-1}\cup\{0\})$.
\end{prop}
\medskip
\medskip

\begin{proof}
 We shall construct such homeomorphism $f$ by induction on $n$.
\medskip

 \textbullet \ The case $n=2$. In this case $B_{2}$ is the unit disk, denoted by $D= \{z\in \mathbb{C}:  |z|\leq 1\}$. Let
$z_{0}= r_{0}e^{i\frac{1}{1-r_{0}^{2}}}\in D$ where $0< r_{0}< 1$ is a real. We shall construct a
$C^{\infty}$-diffeomorphism $f$ on $D$ such that $\omega_{f}(z_{0})= \partial D = S^{1}$
with Fix$(f) =  \partial D\cup \{0\}$. Such a diffeomorphism will be the time-one map of a $C^{\infty}$-flow on $D$. Let $X$ be the vector fields defined on $D$ by
$$X(z)= \begin{cases} e^{\dfrac{i}{1-|z|^{2}}}\left (|z|^{2}+
2i(1-|z|^{2})|z|^{4}\right)\varphi (z), & \textrm{ if } |z|< 1\\
0, & \textrm{ if } |z|=1
\end{cases} $$ where
$$\varphi(z) = \begin{cases} e^{\dfrac{-1}{1-|z|^{2}}},& \textrm{ if } |z|< 1\\
                                0, & \textrm{ if } |z|= 1
                                                 \end{cases}$$

%and %$$F(z)= \dfrac{1}{1-|z|^{2}}, \ \textrm{for } |z|< 1 .$$
%Write $h(r)= \begin{cases} e^{\dfrac{-1}{1-r^{2}}},& \textrm{ if } 0\leq r< 1\\
                               % 0, & \textrm{ if } r=1
                                                 %\end{cases}$. One has $\varphi (z)= h(|z|)$.
 The map $\varphi$ is $C^{\infty}$ on $D$. Indeed, set $h(x)=  e^{\dfrac{-1}{1-x^{2}}}$, $x\in [0,1[$. The derivative of $h$ of order $n$,
($n\in \mathbb{Z}_{+}$)
 is: $h^{(n)}(x)= P_{n}(x)\dfrac{e^{{\dfrac{-1}{1-x^{2}}}}}{(1-x^{2})^{2n}}$, where $P_{n}$ is a polynomial function. Hence
$\underset{x\mapsto 1}\lim h^{(n)}(x) =0$ and so $h$ is $C^{\infty}$ on $1$. We conclude that $\varphi$ is $C^{\infty}$ on $D$ and so on $X$. The set of singular points of $X$ is $\{0\}\cup S^{1}$. Let $(\phi_{t})_{t\in \mathbb{R}}$ denote the flow on $D$ associated to $X$.
Write $\phi_{t}(z_{0}) := z(t)= r(t)e^{i\theta(t)}, \ t\in \mathbb{R}$ where $z(0)= z_{0}$, be the solution of the differential
equation $z^{\prime}(t)= X(z(t))$ with  $z(0)= z_{0}$. Then $r$ is a solution of the real differential equation
$r^{\prime}(t)= r^{2}(t)\varphi(r(t))$ with $r(0)=  r_{0}$ and $\theta(t) = \dfrac{1}{1-r^{2}(t)}$.
 The function $r$ satisfies  $\underset{t\to +\infty}\lim r(t)=1$,
$\underset{t\to -\infty}\lim r(t)=0$. The function $\theta$ is increasing and satisfies $\underset{t\to +\infty}\lim
\theta(t)= +\infty$ and $\underset{|t|\to +\infty}\lim\theta^{\prime}(t)
= \underset{|t|\to +\infty}\lim \dfrac{2r^{3}\varphi(r(t))}{(1-r^{2}(t))^{2}}= 0$.
We let $f= \phi_{1}$ be the time one map of the flow. Then $f$ is a $C^{\infty}$-diffeomorphism of $D$.
Let us show that $\omega_{f}(z_{0})= S^{1}$. One has $f^{n}(z_{0})= \phi_{n}(z_{0})$, so $\omega_{f}(z_{0})\subset S^{1}$.
Conversely, if $z\in S^{1}$, there is a sequence
$(t_{n})_{n}\subset \mathbb{R}$ with $t_{n}\to +\infty$ such that $\phi_{t_{n}}(z_{0})\to z$. As
$\underset{n\to +\infty}\lim r(n)=1$ and $\underset{n\to +\infty}\lim
\theta(n+1)-\theta(n) =0$, it follows that $\underset{n\to +\infty}\lim |f^{[t_{n}]}(z_{0})-\phi_{t_{n}}(z_{0})|=0$ where $[t_{n}]$ is
the integer part of $t_{n}$. Therefore $\underset{n\to +\infty}\lim f^{[t_{n}]}(z_{0})=z$ and so $S^{1}\subset \omega_{f}(z_{0})$.
We conclude that $\omega_{f}(z_{0})= S^{1}$.
\bigskip

\textbullet \ Suppose that for $n\geq 2$, there exists a homeomorphism $f_n$ of $B_{n}$ onto itself satisfying the required properties.
Let $f_{n+1}$ denote the map given by:
$$f_{n+1}(x) = \begin{cases}

 \left(\sqrt{1-x_{n+1}^{2}}~f_n(\frac{1}{\sqrt{1-x_{n+1}^{2}}}(x_1,x_2,\dots,x_n)),x_{n+1}\right), & \textrm{if } x_{n+1}\neq \pm 1\\
(0,\dots, \pm 1),  & \textrm{if } x_{n+1} = \pm 1
\end{cases}$$ where $x=(x_1,x_2,\dots,x_n,x_{n+1})\in B_{n+1}$. It is plain that $f_{n+1}$ is continuous and bijective, hence $f_{n+1}$ is a
homeomorphism of $B_{n+1}$ onto itself. In addition, $S^{n}\cup \{0\}\subset \textrm{Fix}f_{n+1}$. For any $x\in (x_1,x_2,\dots,x_n,x_{n+1})\in B_{n+1}$, we have
$$\omega_{f_{n+1}}(x) = \sqrt{1-x_{n+1}^{2}}~\omega_{f_n}\left (\frac{1}{\sqrt{1-x_{n+1}^{2}}}(x_1,x_2,\dots,x_n)\right)\times\{x_{n+1}\}.$$
We conclude that there is a point $x_{0}\in B_{n+1}\setminus (S^{n}\cup \{0\})$ such that $\omega_{f_{n+1}}(x_{0})\subset S^{n}$ and
 is homeomorphic to the circle $S^1$.
\end{proof}
\medskip

\textit{Proof of Theorem \ref{t4}}. Let $V$ be a free topological $n$-ball in $X$. Its boundary $\partial V$
 is homeomorphic to the sphere $S^{n-1}$. By Proposition \ref{prop51},
there exists a homeomorphism $f: V\to V$ and a point $c\in V$ such that $\partial V \cup \{c\}\subset \textrm{Fix}f$ and for some
point $x\in V\setminus\{c\}$, $\omega_f(x)\subset \partial V$ is homeomorphic to the circle $S^1$. In this way, the map $f$ can be easily extended to a continuous map, say $\widetilde{f}$ defined on the
whole space $X$ by fixing all points in $X\setminus V$. Hence $\widetilde{f}$ is homeomorphism of $X$ onto itself, where
$\omega_{\widetilde{f}}(x)$ is infinite and consisting of fixed points.
We conclude that $X$ does not admit the $\omega$-FTP property. \qed
\bigskip

\textit{Proof of Proposition \ref{p3}}. Let $(x_{1},\dots, x_{n})\in X_{1}\times \dots \times X_{n}$ and assume that
$\omega_{f}((x_{1},\dots, x_{n}))\subset P(f)$. Let $i\in\{1,\dots,n\}$ and  $a_{i}\in \omega_{f_{i}}(x_{i})$, then there exists a sequence
$(n_{k}^{(i)})_{k}\to +\infty$ such that
$f_{i}^{n_{k}^{(i)}}(x_{i})\to a_{i}$. Now the sequence  $(f^{n_{k}^{(i)}}((x_{1},\dots, x_{n}))_{k} = (f_{1}^{n_{k}^{(i)}}(x_{1}), \dots,
f_{n}^{n_{k}^{(i)}}(x_{n}))_{k}$ has a subsequence that converges to a point $(b_{1}^{(i)}, \dots, a_{i},\dots, b_{n}^{(i)})\in P(f)$.
 As $P(f)\subset P(f_{1})\times \dots \times P(f_{n})$ then $a_{i}\in P(f_{i})$. We conclude that $\omega_{f_{i}}(x_{i})\subset P(f_{i})$.
By Theorem \ref{thm2}, the $\omega_{f_{i}}(x_{i})$ ($1\leq i\leq n$) are
finite. As $\omega_{f}(x_{1},\dots, x_{n})\subset \omega_{f_{1}}(x_{1})\times \dots \times \omega_{f_{n}}(x_{n})$,
hence $\omega_{f}(x_{1},\dots, x_{n})$ is finite. \qed
\bigskip

Recall that a Peano continuum can be defined as a continuous image of $[0,1]$. In view of Theorem \ref{thm2},
 one may be tempted to ask the following question.
\medskip

\textbf{Question.} Is it true that a Peano continuum has the
$\omega$-FTP property if and only if it is completely regular ?

\bibliographystyle{amsplain}

\bigskip

\begin{thebibliography}{9}
 \bibitem{ag} S.J. Agronsky, J.G. Ceder,  \emph{ What sets can be $\omega$-limit sets in $E^{n}$},
Real Anal. Exchange, \textbf{17} (1991/92), 97--109.
\bibitem{bdr} F. Balibrea, T. Downarowicz, R. Hric, L. Snoha and V. Spitalsky, \emph{ Almost totally disconnected minimal systems}.
Ergodic Theory Dynam. Systems, \textbf{29},  (2009), 737--766.
\bibitem{ba} F. Balibrea, J.L. Garc\'{\i}a Guirao, J.I. Mu\~{n}oz Casado, \emph{ Description of $\omega$-limit sets of a triangular map on
$I^{2}$}, Far East J. Dyn. Syst. \textbf{3 }(2001), 87--101.
%\bibitem{Bal} F. Balibrea, R. Hric, L.  Snoha, \emph{Minimal sets on graphs and dendrites}, Int. J. Bifurcation and Chaos \textbf{13}
% (2003), 1721--1725.
\bibitem{lsB} L.S. Block, W.A. Coppel, Dynamics in  One Dimension, Lecture Notes in Math, 1513. Springer-Verlag, 1992.
\bibitem{hr}  R. Hric, M. M\'{a}lek, \emph{Omega limit sets and distributional chaos on graphs}, Topology Appl. \textbf{153} (14) (2006), 2469--2475.
\bibitem{jl} V. Jim\'{e}nez, J. Sm\'{\i}tal, \emph{$\omega$-limit sets for triangular mappings}, Fund. Math. \textbf{167} (2001), 1--15.
\bibitem{ku} K. Kuratowski, Topology, vol. II, Academic Press, New York, 1968.
 \bibitem{K} P. Kurka,  Topological and symbolic Dynamics, Soc. Math. de France, (2003).
%\bibitem{Ma} J.H. Mai, E.H. Shi, \emph{$\overline{R}=\overline{P}$ for maps of dendrites $X$ with
%$\textrm{Card}(\textrm{End}(X))<c$}, Int. J. Bifurcation and Chaos, \textbf{19} (2009), 1391-1396.
\bibitem{Li} S. Li, \emph{$\omega$-chaos and topological entropy}, Trans. Amer. Math. Soc. \textbf{339} (1993), 243--249.
\bibitem{Nadler} S. B. Nadler, Continuum Theory: An Introduction, (Monographs and Textbooks in Pure and Applied Mathematics, 158). Marcel Dekker, Inc., New York, 1992.
\bibitem{Nagh1} I. Naghmouchi, \emph{Dynamics of monotone graph, dendrite and dendroid maps},
 Internat. J. Bifur. Chaos Appl. Sci. Engrg., \textbf{21}, (2011), 1--11.
\bibitem{sp} V. Spitalsky, \emph{ Omega-limit sets in hereditarily locally connected continua}.
Topology Appl. \textbf{155 }(2008), 1237--1255.
\end{thebibliography}

\end{document}